\documentclass[12pt]{article}
\usepackage{graphicx}
\usepackage{amsmath}
\usepackage{epsf}
\usepackage{graphpap}

\numberwithin{equation}{section}

\def\be{\begin{equation}}
\def\ee{\end{equation}}

\begin{document}

\begin{center}
\Large{Conditional linearizability criteria for a system of third-order ordinary differential equations}

\bigskip
\end{center}
\begin{center}
{\bf F. M. Mahomed$ ^{a}$, I. Naeem$^{a}$ and Asghar Qadir$^{b}$ }
\end{center}
$^a$Centre for Differential Equations, Continuum Mechanics and Applications\\
School of Computational and Applied Mathematics\\
University of the Witwatersrand\\
Wits 2050, South Africa\\
E-mail: Fazal.Mahomed@wits.ac.za

$^b$Centre for Advanced Mathematics and Physics\\
National University of Sciences and Technology\\
Campus of the College of Electrical and Mechanical Engineering\\
Peshawar Road, Rawalpindi, Pakistan and\\
Department of Mathematical Sciences\\
King Fahd University of Petroleum and Minerals\\
Dhahran 31261, Saudi Arabia\\
E-mail: aqadirs@comsats.net.pk

{\bf {Abstract}}\\
We provide linearizability criteria for a class of systems of
third-order ordinary differential equations (ODEs) that is cubically
semi-linear in the first derivative, by differentiating a system of
second-order quadratically semi-linear ODEs and using the original
system to replace the second derivative. The procedure developed
splits into two cases, those where the coefficients are constant and
those where they are variables. Both cases are discussed and
examples given. \\

\section{Introduction}
Linearization plays an important role in the study of nonlinear
differential equations as there is no standard method to solve
nonlinear ODEs. In 1883, Lie \cite{lie1} gave the linearizability
criteria for scalar second-order ordinary differential equations and
he constructed the most general form of second-order ODEs reducible
to linear ones by changing the dependent and independent variables.
Lie also developed point transformation groups \cite{lie2}. To the
best of our knowledge Lie did not provide linearizability criteria
either for second-order systems or for third-order scalar ODEs,
leave alone systems of third order ODEs. Nor does this aspect appear
in the literature. One can use the approximation (perturbation
method) to convert nonlinear ODEs into linear ones \cite{nay} but
one faces problems in the convergence of the series. A similar
problem arises in numerical schemes. Contact transformations were
used by Chern \cite{che1, che2} to reduce scalar third-order
equations into the linear forms $u'''=0$ and $u''+u=0.$ Grebot
\cite{gre} used fibre preserving point transformations for mapping
third-order ODEs to linear form. A generalization of this work was
later performed by Neut and Petitot \cite{neu}. Ibragimov and
Meleshko \cite{ibr} and Maleshko \cite{mal} also studied the
reduction of scalar third-order ODEs to linear forms. Whereas Neut
and Petitot and Ibragimov and Meleshko used the classical Lie
method, Meleshko dealt with a sub-class of those equations in which
the independent variable did not appear explicitly and hence one
could convert to treating the function to be determined as the
independent variable and and its derivative as the dependent
variable, thereby reducing the order by one.

The linearizability criteria for a system of second order
quadratically semi-linear ODEs were discussed by two of us
\cite{mah}, using the connection between the symmetries and the
isometries of the system of the geodesic equations \cite{tob}. The
criteria require that the curvature tensor be zero, i.e. the space
be flat. Subsequently \cite{qad} a projection procedure of Aminova
and Aminov \cite{ami} was followed, using the translational
invariance of the geodetic parameter in the geodesic equation. A
system of two dimensions was used to get a single cubically
semi-linear ODE and Lie's results on linearization of scalar second
order ODEs was also re-derived. Moreover, invariant criteria for a
system of two cubically semi-linear second-order ODEs to be
reducible to a linear system were obtained.

A method was developed \cite{maq} The purpose of this paper is to
construct the linearizability conditions for a system of third-order
ODEs, cubically semi-linear in the first derivative by
differentiating the system of second order quadratically semi-linear
ODEs and then setting the original system to replace the second
derivative. The procedure employed here uses the derivative of the
vector of the coefficients in the equation and solves a system of
non-homogeneous algebraic equations for it. If this vector becomes
zero, the procedure of solving the system by inverting a matrix
cannot be used directly. As such, one has to go back to the system
and solve it directly. This allows us to construct the
linearizability conditions in all cases. Since there is no work on
explicit linearizability criteria of systems of third order
non-linear ODEs, we cannot compare our results with other works.
However, there is a result for the number of classes of linearizable
systems of ODEs, which is based on there being $3n$ arbitrary
constants for the system, using the classical Lie method. The class
of linearizable equations obtained here is not contained in those
classes as it is non-classical and does not have the required number
of arbitrary constants.

The outline of the work is as follows. In the next section we
briefly mention the geometrical notation used and review the
linearizability criteria for second order quadratically semi-linear
systems. In section 3 we discuss the linearizability criteria for
the general case for a system of third-order ODEs when at least one
of the coefficients is variable. The special case when all the
coefficients are constant is discussed in the next section. In
section 5 we illustrate via some examples to verify our criteria.
Concluding remarks are given in section 6.

\section{Preliminaries}

We shall use the Einstein summation convention that repeated indices
are summed over. The position vector, $x^j$, $(j=1,...,n)$ is
assumed to be a function of a geodetic parameter $s$ and $'$
represents the derivative relative to $s$. The metric tensor will be
denoted by $g_{ij}(x^m)$ and the partial derivative relative to
$x^j$ by ",j". The Christoffel symbols $\Gamma ^i_{jk}$ are defined
in terms of the metric tensor by
\begin{equation}
\Gamma^i_{jk}=\frac{1}{2}g^{il}(g_{jl,k}+g_{kl,j}-g_{jk,l}).
\end{equation}
The geodesic is given by the equation
\begin{equation}
x^{a''}+\Gamma ^a_{bc}x^{b'}x^{c'}=0.
\end{equation}

The Riemann tensor is defined by
\begin{equation}
 R^i_{jkl}=\Gamma ^i_{jl,k}-\Gamma ^i_{jk,l}+\Gamma ^i_{mk}\Gamma ^m_{jl}
 - \Gamma ^i_{ml}\Gamma ^m_{jk},
\end{equation}
which satisfies
\begin{equation}
 R^i_{jkl}=-R^i_{jlk},
 \end{equation}
\begin{equation}
 R^i_{jkl}+R^i_{klj}+R^i_{ljk}=0,
 \end{equation}
and
\begin{equation}
 R^i_{jkl;m}+R^i_{jlm;k}+R^i_{jmk;l}=0.
 \end{equation}
The Riemann tensor in complete covariant form is
\begin{equation}
 R_{ijkl}= g_{im}R^m_{jlk},
 \end{equation}
and has the additional property
\begin{equation}
 R_{ijkl}= -R_{jikl}.
 \end{equation}

A general system of quadratically semi-linear second order ODEs can
be written as
\begin{equation}
x^{a''}+\gamma ^a_{bc}x^{b'}x^{c'}+ \beta^a_b x^{b'}+\alpha^a=0,
\end{equation}
where $\gamma ^a_{bc}, \beta^a_b, \alpha^a$ are functions of the $n$
independent variables. If $\beta^a_b, \alpha^a=0$ we will call it of
{\it geodesic type} if the coefficients $\gamma ^a_{bc}$ can be
written as $\Gamma ^a_{bc}$ for some metric tensor $g_{ab}$.

For example, a two-dimensional system of geodesic type is of the
form
\begin{equation}
\begin{cases}
x''=a(x,y)x'^2+2b(x,y)x'y'+c(x,y)y'^2,\\
y''=d(x,y)x'^2+2e(x,y)x'y'+f(x,y)y'^2,
\end{cases}
\end{equation}
where we have written $x^1=x$, $x^2=y$ and
\begin{equation}
\Gamma ^1_{11}=-a,\quad \Gamma ^1_{12}=-b,\quad\Gamma ^1_{22}=-c,
\end{equation}
\begin{equation}
\quad \Gamma ^2_{11}=-d,\quad \Gamma ^2_{12}=-e,\quad\Gamma
^2_{22}=-f.
\end{equation}

The linearizability condition of \cite{mah} is that the curvature
tensor constructed from these coefficients, regarded as Christoffel
symbols, be zero, i. e. $R^i_{jkl}=0$. In terms of the coefficients
of (3.2) we have
\begin{equation} a_y-b_x+be-cd=0,\end{equation}
\begin{equation} b_y-c_x+(ac-b^2)+(bf-ce)=0,\end{equation}
\begin{equation} d_y-e_x-(ae-bd)-(df-e^2)=0,\end{equation}
\begin{equation} (b+f)_x=(a+e)_y. \end{equation}

\section{Linearization for third-order ODEs}

To obtain the criteria for the third order semi-linear system of
ODEs we differentiate (2.2) to get
\begin{equation}
 x^{a'''}+2\Gamma ^a_{bc}x^{b''}x^{c'}
 +\Gamma ^a_{bc,d}x^{b'}x^{c'}x^{d'}=0.
\end{equation}
This is a total derivative system and may seem to be an artificial
way of getting the third order systems. However, if we now replace
the $x^{b''}$ in (3.1) using (2.2) we get the system
\begin{equation}
{x^a}^{'''}+[\Gamma^a_{(bc,d)}-2\Gamma^a_{p(b}\Gamma^p_{cd)}]
{x^b}^{'}{x^c}^{'}{x^d}^{'}=0,
\end{equation}
which is {\it not} a total derivative system. Thus the general form
will be a third order system of ODEs that is cubically semi-linear
in the first derivative
\begin{equation}
x^{a'''}+A^a_{bcd}x^{b'}x^{c'}x^{d'}=0,
\end{equation}
which will be linearizable if we can make the identification
\begin{equation}
A^a_{bcd}=\Gamma^a_{(bc,d)}-2\Gamma^a_{p(b}\Gamma^p_{cd)}.
\end{equation}

In the case of a system of two third-order ODEs the linearizability
conditions (2.13) - (2.16) are satisfied. Using (3.3) and (3.4) with
the help of (2.13) - (2.16), we obtain
\begin{equation}
 A^1_{111}=-a_x-2a^2-2bd,
 \end{equation}
\begin{equation}
 A^1_{112}=-a_y-2b_x-2(3ab+2be+cd),
 \end{equation}
\begin{equation}
 A^1_{122}=-(2b_y+c_x)-2(ac+2b^2+bf+2ce),
 \end{equation}
\begin{equation}
 A^1_{222}=-c_y-2(bc+cf),
 \end{equation}
\begin{equation} A^2_{111}=-d_x-2ad-2de,
\end{equation}
\begin{equation}
 A^2_{112}=-(d_y+2e_x)-2(2bd+ae+2e^2+df),
 \end{equation}
\begin{equation}
 A^2_{122}=-(2e_y+f_x)-2(cd+2be+3ef),
 \end{equation}
\begin{equation}
 A^2_{222}=-f_y-2ce-2f^2.
 \end{equation}
Invoking (2.13) - (2.16) we can write (3.5) - (3.12) as
\begin{eqnarray}
a_x&=&-(A^1_{111}+2a^2+2bd),\;a_y=-(A^1_{122}/3+2ab+2be), \nonumber \\
b_x&=&-(A^1_{112}/3+cd+2ab+be),\; b_y=-(A^1_{122}/3+b^2+ac+ce+bf), \nonumber \\
c_x&=&-(A^1_{122}/3+2b^2+2ce),\; c_y=-(A^1_{222}/3+2bc+2cf), \nonumber \\
d_x&=&-(A^2_{111}+2ad+2de),\; d_y=-(A^2_{112}/3+2bd+2e^2), \nonumber \\
e_x&=&-(A^2_{112}/3+df+bd+ae+e^2),\; e_y=-(A^2_{122}/3+cd+be+2ef), \nonumber \\
f_x&=&-(A^2_{122}/3+2be+2ef),\; f_y=-(A^2_{222}+2ce+2f^2).
\label{(2)}
\end{eqnarray}
Using the compatibility conditions and for convenience writing
\[A^1_{111}=P,\quad A^1_{112}=Q,\quad A^1_{122}=R, \quad A^1_{222}=S ,\]
\begin{equation}
 A ^2_{111}=T,\quad A ^2_{112}=U,\quad A^2_{122}=V,\quad A^2_{222}=W,
\end{equation}
the system (\ref{(2)}) reduces to
\begin{eqnarray}
 3P_y-2aQ-2dR=Q_x-6bP-2eQ,  \nonumber \\
 Q_y-aR-3dS=R_x-3cP-fQ, \nonumber \\
 R_y-2bR-6eS=3S_x-2cQ-2fR,  \nonumber \\
 3T_y-2aU-2dV=U_x-6bT-2eU,  \nonumber \\
 U_y-aV-3dW=V_x-3cT-fU,  \nonumber \\
 V_y-2bV-6eW=3W_x-2cU-2fV. \label{(3)}
\end{eqnarray}

For $n=2$ the system (3.4) is linearizable if it satisfies
conditions (3.13) and (3.15), where $a,\;b,\;c,\;d,\;e$ and $f$ are
given in the Appendix. Our procedure puts the derivatives of the
$A^i_{jkl}$ as a vector on the right side of a system of linear
equations and solves for the constraints on the coefficients by
inverting the matrix of coefficients on the left side. If the vector
on the right is zero, we cannot use this procedure to obtain the
conditions on the coefficients. As such, for this procedure to work
at least one of the coefficients in the equations must be variable.
If all the coefficients in system (3.4) are constants then the
derivatives of $A^i_{jkl}$ are zero. For this case we need to
develop a separate procedure that is done in the next section.

\section{Special cases}

In the case $A^i_{jkl,m}=0$, the system (3.15) cannot be used for
obtaining the coefficients by inverting the matrix. We have to solve
the equations without resorting to that procedure. This system
reduces to
\begin{eqnarray}
-aQ-dR+3bP+eQ=0, -aR-3dS+3cP+fQ=0, -bR-3eS+cQ+fR=0, \nonumber \\
-aU-dV+3bT+eU=0\nonumber, -aV-3dW+3cT+fU=0, -bV-3eW+cU+fV=0,
\label{(3)}
\end{eqnarray}
and equations (3.5) - (3.12) become
\begin{eqnarray}
 A^1_{111}+2(a^2+bd)=0,\;A^1_{112}+6(ab+be)=0\nonumber,\\
 A^1_{112}+3(cd+2ab+be)=0,\;A^1_{122}+3(b^2+ac+ce+bf)=0\nonumber,\\
 A^1_{122}+6(b^2+ce)=0,\;A^1_{222}+2(bc+cf)=0\nonumber,\\
 A^2_{111}+2(ad+de)=0,\; A^2_{112}+6(bd+e^2)=0\nonumber,\\
 A^2_{112}+3(df+bd+ae+e^2)=0,\;A^2_{122}+3(cd+be+2ef)=0\nonumber,\\
 A^2_{122}+6(be+ef)=0,\;A^2_{222}+2(ce+f^2)=0.\label{(4)}
 \end{eqnarray}

Non-trivial solutions of (4.1) exist if
$$\begin{vmatrix}
-2Q &  6P  &  0 & -2R & 2Q  & 0  \\
-R  &  0   & 3P & -3S & 0   & Q  \\
0   & -2R  & 2Q &  0  & -6S & 2R \\
-2U &  6T  & 0  & -2V & 2U  & 0  \\
-V  &  0   & 3T & -3W & 0   & U  \\
0   & -2V & 2U &  0  & -6W & 2V  \\
\end{vmatrix}=0. $$
This gives the linearizability criteria of (3.4) when all the
coefficients are constants.

From system (\ref{(4)}), we have
\begin{eqnarray}
cd-be=0,\;b^2+ce-ac-bf=0,\;bd+e^2-df-ae=0, \nonumber\\
A^1_{111}=-2(a^2+bd)\nonumber,\;A^1_{112}=-6(ab+be)\nonumber,\\
A^1_{122}=-6(b^2+ce)\nonumber,\;A^1_{222}=-2(bc+cf)\nonumber,\\
A^2_{111}=-2(ad+de)\nonumber,\;A^2_{112}=-6(bd+e^2)\nonumber,\\
A^2_{122}=-6(be+ef),\;A^2_{222}=-2(ce+f^2).\label{(5)}
\end{eqnarray}

In order to solve system (\ref{(5)}), the following cases are
considered.

{\itshape{Case 1}}: $b=0$.\\
The following subcases arise.\\
{\itshape{Case 1.1}}: $c=0,\;d=0.$\\
{\itshape{Case 1.1.1}}: $e=0.$\\
In this case we have
\begin{equation}
A^1_{112}=0,\;A^1_{122}=0,\;A^1_{222}=0,
\;A^2_{111}=0,\;A^2_{112}=0,\;A^2_{122}=0,
\end{equation}
\begin{equation}
\;A^1_{111}=-2 a^2\;A^2_{222}=-2f^2.
\end{equation}
Hence all the $A^i_{jkl}$ are zero except $A^1_{111}$ and
$A^2_{222}$, which must be non-positive since $a$ and $f$ are real. \\
{\itshape{Case 1.1.2}}: $e\not=0.$\\
Simple manipulations yield
\begin{equation}
A^1_{112}=
0,\;A^1_{122}=0,\;A^1_{222}=0,\;A^2_{111}=0,
\end{equation}
\begin{equation}
A^1_{111}=-2
a^2,\;A^2_{112}=-6a^2,\;A^2_{122}=-6af,\;A^2_{222}=-2f^2,
\end{equation}\\
which yield\\
\[a^2=-A^1_{111}/2,\;f^2=-A^2_{222}/2, \]
\begin{equation}
{(A^2_{122})}^2-9 A^1_{111}A^2_{222}=0,\;A^2_{112}-3A^1_{111}=0.
\end{equation}
This further requires $A^1_{111},\;A^2_{112},\;A^2_{222}\leq 0$.

{\itshape{Case 1.2}}: $c=0,\;d\not=0.$\\
Straightforward calculations lead to
\begin{equation}
e^2-df-ae=0,
\end{equation}
and
\begin{equation}
A^1_{112}= 0,\;A^1_{122}=0,\;A^1_{222}=0,
\end{equation}
\begin{equation}
A^1_{111}=-2 a^2,\;A^2_{111}=-2(ad+de),
\end{equation}
\begin{equation}
A^2_{112}=-6e^2,\;A^2_{122}=-6ef,\;A^2_{222}=-2f^2,
\end{equation}\\
which yield
\[a^2=-\frac{A^1_{111}}{2},\;d^2=\frac{{3(A^2_{111}})^2}{2[-3A^1_{111}-A^2_{112}\pm
2\sqrt{3A^1_{111}A^2_{112}}]} ,\]
\begin{equation}
e^2=-A^2_{112}/6,\;f^2=-A^2_{222}/2,
\end{equation}
with conditions
\begin{equation}
{(A^2_{122})}^2-3A^2_{112} A^2_{222}=0,
\end{equation}
\begin{equation}
3A^2_{111}=[\pm \sqrt{\frac{A^1_{111}}{A^2_{222}}}+\frac{A^2_{112}}
{A^2_{122}}][A^2_{112}\pm \sqrt{3A^1_{111}A^2_{112}}].
\end{equation}
Thus we require that $A^1_{111},\;A^2_{112}\leq 0,\;A^2_{222}< 0$.

{\itshape{Case 1.3}}: $c\not=0,\;d=0.$ \\
We find that $e=a$ and
\begin{equation}
A^1_{112}= 0,\;A^2_{111}=0,
\end{equation}
\begin{equation}
A^1_{111}=-2 a^2,\;A^1_{122}=-6ac,\;A^1_{222}=-2cf,
\end{equation}
\begin{equation}
A^2_{112}=-6a^2,\;A^2_{122}=-6af,\;A^2_{222}=-2(ac+f^2),
\end{equation}\\
which yield
\begin{equation}
a^2=-\frac{A^1_{111}}{2},\;c^2=-\frac{{(A^1_{122}})^2}{18A^1_{111}},\;
f^2=-\frac{{(A^2_{122}})^2}{18A^1_{111}},
\end{equation}
\begin{equation}
A^2_{112}-3A^1_{111}=0,\;9A^1_{111}A^2_{222}-3A^1_{111}A^1_{122}-{(A^2_{122})}^2=0,
\end{equation}
\begin{equation}
81{(A^1_{111})}^2{(A^1_{222})}^2-{(A^1_{122})}^2{(A^2_{122})}^2=0.
\end{equation}
In this case we need to have $A^1_{111}<0,\;A^2_{112}\leq 0.$

{\itshape{Case 2}}: $b\not=0$.\\
The subcases are as follows.

{\itshape{Case 2.1}}: $d=0.$\\
In this case we find that
\begin{equation}
e=0,\;b^2-ac-bf=0,
\end{equation}
and
\begin{equation}
A^2_{111}=0,\;A^2_{112}=0,\;A^2_{122}=0,
\end{equation}
\begin{equation}
A^1_{111}=-2 a^2,\;A^1_{112}= -6ab,\;A^1_{122}=-6b^2,
\end{equation}
\begin{equation}
A^1_{222}=-2(bc+cf),\;\;A^2_{222}=-2f^2,
\end{equation}\\
which imply that
\begin{eqnarray}
a^2=-A^1_{111}/2,\;b^2=-A^1_{122}/6,\;f^2=-A^2_{222}/2, \nonumber\\
c^2=\frac{{3(A^1_{222}})^2}{2[-3A^2_{222}-A^1_{122}\pm
2\sqrt{3A^1_{122}A^2_{222}}]}
\end{eqnarray}
with extra conditions on $A^i_{jkl}$ given by
\begin{equation}
{(A^1_{112})}^2-3A^1_{111}A^1_{122}=0,
\end{equation}
\begin{equation}
3A^1_{222}=[\pm \sqrt{\frac{A^2_{222}}{A^1_{111}}}+\frac{A^1_{122}}
{A^1_{112}}][A^1_{122}\pm \sqrt{3A^1_{122}A^2_{222}}].
\end{equation}
and $A^1_{111}<0,\;A^1_{122},\;A^2_{222}\leq 0.$

{\itshape{Case 2.2}}: $d\not=0,\; e=cd/b.$\\
We find that
\begin{equation}
b^3+c^2d-b^2f-abc=0,\; b^3d+c^2d^2-b^2df-abcd=0,
\end{equation}
and
\begin{equation}
A^1_{111}=-2(a^2+bd),\;A^1_{112}=-6(ab+be),\;A^1_{122}=-6(b^2+ce),
\end{equation}
\begin{equation}
A^1_{222}=-2(bc+cf),\;A^2_{111}=-2(ad+de),\;A^2_{112}=-6(bd+e^2),
\end{equation}
\begin{equation}
A^2_{122}=-6(be+ef),\;A^2_{222}=-2(ce+f^2).
\end{equation}

In order to construct the linearizability conditions for above
system one has to replace $e=\lambda \frac{d}{b}$ or $e=\lambda
\frac{c}{b}$ to get a $5\times5$ matrix and the determinant of the
resultant matrix must be zero.

\section{Examples}
To test our linearization criteria, we utilize the following examples.

{\bf 1.} For the case when at least one of the coefficients is
variable e.g. $P=0,\;Q=0,\;R=3,\;S=0,\;T=0,\;U=-6/x^2,\;V=0,\;W=2,$
the system of two third-order equations is linearizable as it
satisfies conditions (3.13)-(3.15) with $a,b,c,d,e$ and $f$ given in
the Appendix.\\

The following set of six equations can be obtained from \cite{mah}
by using the values of $a,b,c,d,e$ and $f$
\begin{equation}
p_x=0,\;q_x=q/x,\;r_x=2r/x,\nonumber
\end{equation}
\begin{equation}
p_y=2q/x,\;q_y=-xp+r/x,\;r_y=-2xq.
\end{equation}
A solution of system (5.1) is
\begin{equation}
p=1,\;q=0,\;r=x^2.\nonumber
\end{equation}
To obtain the Cartesian coordinates $u(x,y)$ and $v(x,y)$, one has
to solve (\cite{mah})
\begin{equation}
u_x^2+v_x^2=1,\nonumber
\end{equation}
\begin{equation}
u_xu_y+v_xv_y=0,\nonumber
\end{equation}
\begin{equation}
u_y^2+v_y^2=x^2,
\end{equation}
Thus
\begin{equation}
u=x\cos y,\nonumber
\end{equation}
\begin{equation}
v=x\sin y,\nonumber
\end{equation}
are the transformations that linearize the system.

{\bf 2.} The system which has
$P=0,\;Q=-3,\;R=0,\;S=0,\;T=0,\;U=9y^2,\;V=18,\;W=-6/y^2$, is
linearizable since it holds conditions (3.13)-(3.15) with
$a,\;b,\;c,\;d,\;e$ and $f$ given in the Appendix. In the same
manner as in example 1, one can obtain
\begin{equation}
p=1+x^2-2x/y+1/y^2,\nonumber
\end{equation}
\begin{equation}
q=(1+x^2)/y^2-x/y^3,\nonumber
\end{equation}
\begin{equation}
r=(1+x^2)/y^4.\nonumber
\end{equation}
For the coordinates transformations one can solve system (5.2) with
the help of $p,\;q,\;r$ to obtain
\begin{equation}
u=x-1/y,\nonumber
\end{equation}
\begin{equation}
v=x^2/2-x/y.\nonumber
\end{equation}

{\bf 3.} The system of two third order equations with constant
coefficients $P=-2,\;Q=0,\;R=-6,\;S=0,\;T=0,\;U=-6,\;V=0,\;W=-2$ is
also linearizable as the determinant of (\ref{(3)}) is zero.\\
The coordinate transformation for the above system which gives the
linearization is
\begin{equation}
u=c_1 e^{-y-x}+c_2 e^{y-x},
\end{equation}
\begin{equation}
v=c_1e^{-y-x}-c_2e^{y-x}.
\end{equation}
Note that there are four arbitrary constants which will appear in
the solution for each of the above examples and therefore we have
conditional linearizability subject to a system of two second-order
equations.

\section{Concluding Remarks}

In this paper we have provided conditional linearizability criteria
for a class of third-order systems of ODEs. To the best of our
knowledge the linearization of systems of third-order ODEs has not
been studied before in the literature. The system of two third-order
ODEs that is cubically semi-linear in the first derivative was
obtained by differentiating a system of two quadratically
semi-linear ODEs and then we used the original system to replace the
second derivative. The resultant system was then not a total
derivative.

The criteria developed are discussed in two parts. In the first part
we constructed the linearizability conditions when at least one of
the coefficients is a non-constant function, so that the system
obtained for the linearizability conditions be invertible. For the
constant coefficients case we obtain a system which is independent
of the derivatives. The non-trivial solution of the resulting system
requires that the determinant of the coefficients be zero. There
were various sub-cases. The linearizability conditions for each of
the sub-cases of constant coefficients were constructed. This
procedure will hopefully give rise to further studies in the
construction of linearizability criteria for systems of higher-order
equations.

{\bf Acknowledgements}\\
AQ is most grateful to DECMA and the School of Computational and
Applied Mathematics, University of the Witwatersrand and for some
useful comments by Profs. P. Leach, S. Maleshko and R. Popovych.

{\bf \Large{Appendix}}\\\\
By using Mathematica one can find $a,\;b,\;c,\;d,\;e$ and $f$ from
system (3.15) in terms of $P,\;Q,\;R,\;S,\;T,\;U,\;V$ and $W$ as
{\[a=\frac{1}{\Delta }[3(-9 S^2TU^2 + RSU^3 + 27S^2 T^2V -QSU^2V -
3QSTV^2 + 3PSUV^2\]
 \[- 27RST^2 W + 18 Q S T U W - Q R U^2 W + 3 Q R T V W - 27P S T V W + Q^2 U V W\]
\[- 3 PR U V W -9 Q^2 T W^2 + 27 P R TW^2)(3 P_y - Q_x)) + 2 (9 R S T U^2 - R^2 U^3\]
\[ - 27 R S T^2 V + 2 Q R U^2 V - 9 P S U^2 V + 27 P S T V^2 - Q^2 U V^2 + 27 R^2 T^2 W\]
\[ - 18 Q R T U W + 9 P R U^2 W + 9 Q^2T V W - 27 P R T V W)( Q_y - R_x)\]
\[ + 3(Q T - P U)(-3S U^2 + 9 S T V + R U V - Q V^2 - 9 R T W + 3 Q U W)( R_y - 3  S_x)\]
\[ - 3 ((QR S U^2 - 9 P S^2 U^2 - 3 Q R S T V + 27 P S^2 T V - Q^2 S U V + 3 P R S U V \]
\[ + 3 Q R^2 T W- 27 P R S T W - Q^2 R U W - 3 P R^2 U W + 18 P Q S U W + Q^3 V W \]
\[ - 27 P^2 S V W - 9 P Q^2 W^2+ 27 P^2 R W^2)(3 T_y - U_x)+ 2 (Q R^2 U^2 \]
\[- 9 P R S U^2 - 9 Q^2 S T V + 27 P R S T V - 2 Q^2 R U V + 18 P Q SU V + Q^3V^2\]
\[ - 27 P^2 S V^2 + 9 Q^2 R T W - 27P R^2 T W - 9 P Q^2 V W + 27 P^2 R V W)(U_y -V_x)\]
$$
- 3 (Q T - P U)(R^2 U - 3 Q S U - Q R V + 9 P SV + 3 Q^2 W - 9 P R
W)( V_y - 3W_x)],\eqno{(A.1)}
$$\\
\[b=\frac{1}{\Delta }[-3(3S^2U^3-9S^2TUV-2RSU^2V+3RSTV^2+2QSUV^2-3PSV^3+9RSTUW\]
\[+R^2U^2W-6QSU^2W-3R^2TVW+9PSUVW-Q^2V^2W+3PRV^2W+3Q^2UW^2\]
\[-9PRUW^2)(3P_y-Q_x)-2(R U - QV)(-3SU^2+9STV+RUV-QV^2-9RTW\]
\[+3QUW)(Q_y  - R_x) -(RU-QV)(RU^2-3RTV-QUV+3PV^2+9QTW\]
\[-9PUW)(R_y-3S_x)-3(R^2SU^2-3QS^2U^2-3R^2STV+9QS^2TV-Q^2SV^2\]
\[+3PRSV^2+3R^3TW-9QRSTW-2QR^2UW+6Q^2SUW+2Q^2RVW-3PR^2VW\]
\[-9PQSVW-3Q^3W^2+9PQRW^2)(3 T_y- U_x)+2(RU-QV)(R^2U-3QSU-QRV\]
\[+9PSV+3Q^2W-9PRW)(U_y-V_x)-(RU-QV)(3R^2T-9QST-QRU\]
$$
+9PSU+Q^2V-3PRV) (V_y-3W_x)],\eqno{(A.2)}
$$\\
\[ c=\frac{1}{\Delta }[-3(SV-RW)(3SU^2-9STV-RUV+QV^2+9RTW-3QUW)(3P_y-Q_x)\]
\[-6(-SV+RW)(RU^2-3RTV-QUV+3PV^2+9QTW-9PUW)(Q_y-R_x)\]
\[-(R^2U^2V-9QSTV^2-2QRUV^2+9PSUV^2+Q^2V^3-9R^2TUW+27QSTUW\]
\[-27PSU^2W+18QRTVW-9PQV^2W-27Q^2TW^2+27PQUW^2)(R_y-3S_x)\]
\[+3(-SV+RW)(R^2U-3QSU-QRV+9PSV+3Q^2W-9PRW)(3T_y-U_x)\]
\[-6(3R^2T-9QST-QRU+9PSU+Q^2V -3PRV)(-SV+RW)(U_y-V_x)\]
\[+(-9R^2STU+27QS^2TU+R^3U^2-27PS^2U^2-2QR^2UV+18PRSUV+Q^2RV^2\]
$$
-9PQSV^2+9QR^2TW-27Q^2STW-9PR^2UW+27PQSUW)(V_y-3W_x)],\eqno{(A.3)}
$$\\
\[d=\frac{1}{\Delta }[9RSTU^2- R^2U^3-27RST^2V+2QRU^2V-9PSU^2V+27PSTV^2-Q^2UV^2\]
\[+27R^2T^2W-18QRTUW+9PRU^2W+9Q^2TVW-27PRTVW)(3P_y-Q_x)\]
\[+6(QT-PU)(-3SU^2+9STV+RUV-QV^2-9RTW+3QUW)(Q_y- R_x)\]
\[+3(QT-PU)(RU^2-3RTV-QUV+3PV^2+9QTW-9PUW)(R_y-3S_x)+QR^2U^2\]
\[-9PRSU^2-9Q^2STV+27PRSTV-2Q^2RUV+18PQSUV+Q^3V^2-27P^2SV^2\]
\[+9Q^2RTW-27PR^2TW-9PQ^2VW+27P^2RVW)(3T_y-U_x)- 6(QT\]
\[-PU)(R^2U-3QSU-QRV+9PSV+3Q^2W-9PRW)(U_y - V_x)\]
$$
+ 3(QT - PU)(3R^2T-9QST -QRU+9PSU+Q^2V-3PR
V)(V_y-3W_x)],\eqno{(A.4)}
$$\\
\[e=\frac{1}{\Delta }[-(RU-QV)(-3SU^2+9STV+RUV-QV^2-9RTW+3QUW)(3P_y-Q_x)\]
\[-2(RU-QV)(RU^2-3RTV-QUV+3PV^2+9QTW-9PUW)(Q_y-R_x)\]
\[-3(-R^2TU^2+3QSTU^2-3PSU^3+3R^2T^2V-9QST^2V+9PS TU V\]
\[ + 2PRU^2V+ Q^2TV^2-6PRTV^2-2PQUV^2+3P^2V^3-3Q^2TUW\]
\[+3PQU^2W+9PQTVW-9P^2UVW)(R_y-3S_x)+(RU-QV)(R^2U -3QSU\]
\[-QRV+9PSV+3Q^2W-9PRW)(3T_y -U_x)-2(RU-QV)(3R^2T-9QST\]
\[-QRU+9PSU+Q^2V-3PRV)(U_y-V_x)+3(3R^3T^2-9QRST^2-2QR^2TU\]
\[+3Q^2STU+9PRSTU+PR^2U^2-3PQSU^2+2Q^2RTV-6PR^2TV-PQ^2V^2\]
$$
+3P^2RV^2-3Q^3TW+9PQRTW+3PQ^2UW-9P^2RUW)(V_y-3W_x)],\eqno{(A.5)}
$$\\
\[f=\frac{1}{\Delta }[3(SV-RW)(RU^2-3RTV-QUV+3PV^2+9QTW-9PUW)(3P_y-Q_x)\]
\[-2(R^2U^2V-9QSTV^2-2QRUV^2+9PSUV^2+Q^2V^3-9R^2TUW+27QSTUW\]
\[-27PSU^2W+18QRTVW-9PQV^2W-27Q^2TW^2+27PQUW^2)(Q_y-R_x)\]
\[-3(-R^2TUV+3QSTUV-3PSU^2V+QRTV^2+PRUV^2-PQV^3+9R^2T^2W\]
\[-27QST^2W-3QRTUW+27PSTUW+3PRU^2W-18PRTVW+9P^2V^2W\]
\[+27PQTW^2-27P^2UW^2)(R_y-3S_x)-3(3R^2T-9QST-QRU+9PSU+Q^2V\]
\[-3PRV)(-SV+RW)(3T_y-U_x)+2(-9R^2STU+27QS^2TU+R^3U^2-27PS^2U^2\]
\[-2QR^2UV+18PRSUV+Q^2RV^2-9PQSV^2+9QR^2TW-27Q^2STW-9PR^2UW\]
\[+27PQSUW)(U_y-V_x)-3(-9R^2ST^2+27QS^2T^2+R^3TU-27PS^2TU-QR^2TV\]
\[-3Q^2STV+18PRSTV-PR^2UV+3PQSUV+PQRV^2-9P^2SV^2+3Q^2RTW\]
$$
-27PQSTW-3PQRUW+27P^2SUW)(V_y-3W_x)],\eqno{(A.6)}
$$
where
\[\Delta=(2(-R^3(U^3-27T^2W)-9Q^2(2ST+PW)(V^2-3UW)+27PS(SU^3-3STUV\]
\[-PV^3+3PUVW)+Q^3(V^3-27TW^2)+3R^2(3ST(U^2-3TV)+6P(U^2-3TV)W\]
\[+QU(UV-9TW)+27QS(ST(-U^2+3TV)+P(UV^2-2U^2W-3TVW))\]
\[+3R(Q^2V(-UV+9TW)-3Q(ST-PW)(-UV+9TW)+9P(PW(V^2-3UW)\]
$$
+S(-U^2V+2TV^2+3TUW))))).\eqno{(A.7)}
$$

\end{document}